\documentclass[12pt]{article}

\usepackage{amsmath,amsthm}
\usepackage{amssymb}

\usepackage{enumitem}

\usepackage{graphicx}

\usepackage[T1]{fontenc}

\newtheorem{theorem}{Theorem}[section]
\newtheorem{corollary}[theorem]{Corollary}
\newtheorem{lemma}[theorem]{Lemma}
\newtheorem{proposition}[theorem]{Proposition}

\theoremstyle{definition}
\newtheorem{definition}[theorem]{Definition}
\newtheorem{remark}[theorem]{Remark}

\numberwithin{equation}{section}

\newcommand{\C}{\mathbb{C}}
\newcommand{\R}{\mathbb{R}}
\newcommand{\N}{\mathbb{N}}

\begin{document}


\baselineskip=17pt


\title{On some quasianalytic classes of $C^\infty$ 
functions }

\author{Abdelhafed Elkhadiri\\
University Ibn Tofail\\ 
Faculty of Sciences\\
Department of Mathematics\\
Kenitra, Morocco\\
E-mail: elkhadiri.abdelhafed@uit.ac.ma
}

\date{}

\maketitle


\renewcommand{\thefootnote}{}



\renewcommand{\thefootnote}{\arabic{footnote}}
\setcounter{footnote}{0}

\begin{abstract}
    This expository article is devoted to the notion
     of quasianalytic classes and the Borel mapping.
    Althougt quasianalytic classes are well known 
    in analysis since several decades.
 We are interested in certain
properties of Denjoy-Carleman's quasianalytic 
classes, 
such as the non-surjectivity of the Borel 
mapping, the property of monotonicity.
We try to see if it remains true for other 
quasianalytic classes, such as for example, 
the classes of indefinitely differentiable 
functions definable in a polynomially 
bounded o-minimal structures. What motivated this 
is the fact of having shown in a previous article 
the existence of quasianalytic classes 
where  Borel mapping is surjective.

\end{abstract}
\section{Introduction}
Analytic functions on an interval $[a,b]\subset \R, 
\,\, a <b,$ posses the following two equivalent 
properties:
\begin{enumerate}
\item [$\mathcal{B})$]  An analytic function is 
determined on $[a,b]$ as soon as it is known in 
a subinterval of $[a,b]$.
\item  [$\mathcal{DC})$]An analytic function is 
determined on $[a,b]$ by its value and the values 
of its successive derivatives at a point 
$c\in [a,b]$.
\end{enumerate}
This is a direct consequence of the fact that 
analytic functions are developable in a 
Taylor's series in the neighborhood of the point 
and the Taylor's development identifies the function. 
It was thought for a long time that the analytic 
 functions where the only ones which were 
determined by their values and the values of there 
derivatives at a single point of the interval 
$[a,b]$.  It was
Borel who first proved the existence of functions 
belonging to more general classes than that 
of analytic functions, which possessed the property 
that they are determined by their values and the
values of all their successive derivatives at 
a single point $x_0\in [a,b]$, see \cite{Borel}. Hi gave 
them the name {\it {quasianalytic}} functions.\\
For any subvector space $\mathcal{A}$ of the ring of 
$C^\infty$ functions on the  interval $[a,b]$,
$C^\infty([a,b])$, if condition $\mathcal{DC})$ 
is verified by the elements of $\mathcal{A}$
then condition $\mathcal{B})$  is verified.
The reciprocal is not true in general
in the class of $C^\infty$ functions.
Take the function defined in the interval 
$[0,1]$  by 

$$f(x)=\left\lbrace
\begin{array}{ll}
0 & \mbox{if $ x =0$ }\\
\exp{(-\frac{1}{x})} & \mbox{if}\,\, x\neq 0.
\end{array}
\right.
$$
In other words in the sitting of 
 $C^\infty$ functions,  condition $\mathcal{DC})$ 
 is stronger than condition $\mathcal{B})$.\\
 Condition $\mathcal{B})$ is the definition of 
 quasi-analyticity that Bernstein had adopted, 
 see \cite{Bern1}, while
condition $\mathcal{DC})$ was adopted by 
Denjoy-Carleman. Bernstein's definition of 
quasi-analyticity does not
use the fact that the functions are differentiable,
 and can therefore be defined even for continuous
functions.\\
In this paper we will say that the subvector space
$\mathcal{A} \subset  C^\infty([a,b])$ is 
quasi-analytic, if $\mathcal{A}$ satisfied
condition $\mathcal{DC})$.\\
Let us  mention an other property that will 
interest us in this paper, which we call 
{\it{analytic monotonicity}} property.\\
If $f$  is an analytic functions on the interval
$[a,b]\subset\R$ such that all its derivatives 
are positive at the left end point $a$, then all the 
derivatives are positive in the whole interval.

\section{ Quasi-analytic functions of  
real variable.}
In this section we are concerned with the approach 
used by Denjoy to construct some  of
quasianalytic classes of real functions.
\subsection { Quasi-analytic functions according 
to Denjoy’s point of view.}
Let $f$ be a $C^\infty$ function on the interval 
$[a,b]\subset\R$ and put, for each $n\in\N$,
$$M_n(f) = \sup_{x\in[a,b]}\mid f^{(n)}(x)\mid, 
$$ where 
$f^{(n)}$
is the nth derivative of the function $f$.\\
Recall the characterization, given by Pringsheim, 
of analytic functions on the interval 
$[a,b]\subset\R$. 
Let $f$ be  a $C^\infty$ function on the interval 
$[a,b]\subset\R$, then
$$ f \,\,\,\,\mbox{is analytic in the interval}
\,\,\,\,
[a,b] \Longleftrightarrow 
\sqrt[n]{M_n(f)} \leq C.n,\,\,\,
\forall n\in\N, \eqno{(\mathcal{C}_0)}$$
where $C$ is a positive constant 
independent of the integer 
$n$.\\
Denjoy asked himself if it is possible to 
enlarge the class of analytic functions on $[a,b]$
without losing the condition $\mathcal{DC)}$. 
For this purpose he inserted factors which become 
infinitely great with $n$ to the second member of 
$(\mathcal{C}_0)$ and in this way defined 
 different classes of $C^\infty$ 
functions on
the interval $[a,b]$, ($\mathcal{C}_j)_{j\in \N}$,
 characterized by the following
conditions:
\begin{equation}
    \sqrt[n]{M_n(f)} \leq C.n\log n,  \,\,\,
    \forall n > 1,\,\,\,\, \,\,\,\,\,(\mbox{class} 
    \,\,\,\mathcal{C}_1)
\end{equation}
\begin{equation}
    \sqrt[n]{M_n(f)} \leq C.n\log n .\log(\log n), 
     \,\,\, \forall n > e, ,\,\,\,\, (\mbox{class} 
     \,\,\,\mathcal{C}_2)
\end{equation}
and so on. \\By uniformity we designate by 
$\mathcal{C}_0$  the analytical class, that is to 
say the class of functions characterized by
\begin{equation}
    \sqrt[n]{M_n(f)} \leq C.n,  \,\,\,
    \forall n \in\N,\,\,\,\, (\mbox{class} 
    \,\,\,\mathcal{C}_0)
\end{equation}

He proved that the functions of this different 
classes, $\mathcal{C}_1,\,\, \mathcal{C}_2, ....$ 
were still satisfied the condition $\mathcal{DC})$.
See appendix at the end of this article for an 
introduction and treatment of these quasianalytic
 classes of $C^\infty$ functions
in a  more conveniently way. As we see, 
 we can have almost all 
properties.\\

Denjoy noticed that the reciprocal of the 
second members of the inequalities (2.1),  (2.2) and
(2.3) are the general
terms of divergent series. He was therefore 
led to announce the following theorem without 
prove it.
\begin{theorem}
    Let $f$  be a $C^\infty$  function on the 
    interval $[a,b]$. The function $f$ is completely 
    determined
    in the whole interval $[a,b]$ by its value 
    and the values of its derivatives in any point 
    of $[a,b]$, if the
    series of positive terms:
    $$\sum_{n=0}^\infty \frac{1}{\sqrt[n]{M_n(f)}}$$
    is divergent.

\end{theorem}
In order to show this result, T. Carleman 
considered 
classes of functions more general than those
considered by Denjoy. He proceeds as follows:\\

Let $M = (M_n)_{n\in \N}$  be a sequence of positive
 numbers. 
 We denote by $C_M([a,b])\subset C^\infty ([a,b])$ 
 the class of infinitely differentiable 
 functions in the interval $[a,b]$ satisfying
$$ 
\mid f^{(n)}(x)\mid \leq k. C^n M_n, \,\,\forall 
n\in\N,\,\,\,\forall x\in [a,b],
$$
where $k$ , $C$  are positive constants 
(depending on $f$, but not on $n$).\\
We remark that if $k$ is omitted in the definition,
 then if $n = 0$,  we have
$ \sup\limits_{x\in [a,b]}\mid f(x)\mid \leq M_0$, 
which is restrictive. \\ In the following, 
we suppose $M_0 = 1$. The class $C_M([a,b])$ is a 
vector space. A class that satisfies condition 
$\mathcal{DC})$ will be called quasianalytic class 
in the sense of Denjoy-Carleman.
\begin{proposition}
    Let $ M = (M_n)_{n\in \N}$  be a sequence of 
    positive numbers. The class $C_M([a,b])$ is
    quasianalytic in the sense of Denjoy-Carleman 
    if and only if condition $\mathcal{B})$
     is satisfied.   
\end{proposition}
\begin{proof}
It is enough to show that condition $\mathcal{B})$
 implies that the class $C_M([a,b])$
    is quasianalytic in the sense of 
    Denjoy-Carleman.\\
    Let $f\in C_M([a,b])$ such that $f^{(n)}(c) =0$,
    $\forall n\in\N$, where $c\in ]a,b[$.  Let 
    $g$ be the function defined as follows: 
    $$
    g(x)=\left\lbrace
    \begin{array}{ll}
    0 & \mbox{if } x \in [a,c] \\
    f(x) & \mbox{if}\,\, x\in [c,b].
    \end{array}
    \right.
    $$
    It is a straightforward observation that the 
    function $g$ is a $C^\infty $ function on 
    $[a,b]$ and $g\in C_M([a,b])$.
    By condition
    $\mathcal{B})$, we deduce that the function 
    $g$  is zero and therefore the function 
    $f$  vanishes on the interval $[c,b]$
    and so $f$ is identically zero. \\
    If $c =a$  [resp. $c = b$],  we set the 
    function $g=0$ on $[a -\eta , a]$  
    [resp. $[b,b+\eta]$]
     where $\eta  > 0$  and 
    $g = f$  on $[a, b]$.
\end{proof}
Let’s eliminate some trivially case of the 
    sequence $ M = (M_n)_{n\in \N}$
     where the class  $C_M([a,b])$
     is quasianalytic in the sense of 
     Denjoy-Carleman.
\begin{proposition}
 Let  $ M = (M_n)_{n\in \N}$ be a sequence of 
 positive numbers. Then the class $C_M([a,b])$
 is quasianalytic in the sense of 
 Denjoy-Carleman if $\liminf\limits_{n\to \infty}
 \sqrt[n]{M_n(f)} < \infty$.
   
\end{proposition}
\begin{proof}
    Put $ \lambda = \liminf\limits_{n\to \infty}
    \sqrt[n]{M_n(f)}$.  For every $\epsilon >0$ 
    corresponds an infinite increasing sequence 
    of natural numbers $(n_j)_{j\in\N}$ 
    such that   $M_{n_j} \leq
    (\lambda +\epsilon)^n{_j}$. 
    consider   $c\in [a,b]$ 
     such that $f^{(n)}(c) = 0$, $\forall n\in\N$,
then
$$ \mid f(x)\mid \leq 
\frac{\mid f^{n_j}(\theta)\mid }{n_j !}
\mid x-c\mid ^{n_j} \leq \frac{M_{n_j}}{n_j !}
\mid x-c\mid ^{n_j} \leq \exp(\lambda + \epsilon)
\mid x-c\mid ^{n_j},\forall j\in\N.
$$
where $\theta$  is in the open interval of 
end points  $x$  and $c$.
Suppose $\mid x-c\mid <1$, we see then if
 we let $j$ tend to infinity, this would imply $f=0$
on a subinterval of $[a,b]$ containing the 
point $c$. By doing the
same thing with one end of this interval 
and so on, we show that $f$ is zero on
$[a,b]$. 
\end{proof}
From now on, we assume that
$$
\liminf\limits_{n\to \infty}
    \sqrt[n]{M_n(f)} = \infty.
$$    
Carleman provides a complete answer to the 
Theorem 2.1
\begin{theorem}
    The class $C_M([a,b])$ is quasianalytic if, 
    and only if, $\sum\limits_{n=0}^\infty 
    \frac{1}{\beta_n} =\infty$, where $\beta_n = 
\inf\limits_{k\geq n}\sqrt[k]{M_k(f)}$.
\end{theorem}
\begin{remark}
 Although we are not going to address the 
 question of the reconstitution of a function 
 in a quasianalytic class by means of its 
 derivatives at a given point, let us point 
 out that this question 
 is treated in  \cite[Chapitre VII]{Carleman} and 
 \cite{Carleman1}.
\end{remark}

For the proof of theorem 2.4, it is often convenient
 to deal with other equivalent statements involving 
 some other sequence. \\ 
In order to do that, we need to change the 
sequence $ M = (M_n)_{n\in\N}$ by an other with 
suitable properties, and this change does affect 
the quasianalyticity of the class
$C_M([a,b])$.\\
We introduce here the so-called convex 
regularization by means of the logarithm.
\begin{definition}
    A sequence of positive real numbers 
    $ M = (M_n)_{n\in\N}$ is said to be 
    log-convex, if and only if, for all 
    $n \geq 1$, we have  
    $M_n^2\leq M_{n-1}M_{n+1}$.
\end{definition}
The condition $\liminf\limits_{n\to \infty}
\sqrt[n]{M_n(f)} = \infty$ implies the existence 
of the convex regularization by means
of the logarithm of the sequence $ M = (M_n)_{n\in\N}$,
 see \cite{Man}, that is, a sequence 
 $ M_c = (M_n^c)_{n\in\N}$ such that 
\begin{enumerate}
    \item The function $ k\mapsto \log M^c_k$ 
    is convex, i.e.  $ \log M^c_k \leq \frac{1}{2}
    (\log M^c_{k-1} + \log M^c_{k +1}  )$. 
    \item $M^c_k \leq M_k,\,\,\,\forall k\in\N$
    \item  there is a sequence $ 0 = n_0 < 
    n_1 < n_2 < \ldots < n_k < \ldots $
    , called the {\it{ principal sequence}}, such that
    $M^c_{n_k}  = M_{n_k},\,\,\forall k\in\N$
    and the function $k\mapsto \log M^c_{n_k}$ is 
    linear in each $[n_k, n_{k+1}]$.
\end{enumerate}
The convex regularized sequence by means of 
the logarithm is the largest convex minorant of 
the function $n\mapsto \log M_{n}$. 
We give an idea of the construction of such sequence.

Let's  first recall the definition of 
Newton's  polygon attached to the sequence 
$(\log M_n)_{n\in \N}$  (under
the condition $\liminf\limits_{n\to \infty}
\sqrt[n]{M_n(f)} = \infty$).\\

Consider in the plane $x0y$ the points 
$P_n  = (n, \log M_n)$, $n\in\N$. Let 
the half-line passing through the point
$P_0  = (0, \log M_0)$ and pointing to the 
negative direction of $0y$. Let's turn our 
half-line in sense counter clock wise until 
it meets a point $P_n  = (n, \log M_n)$. 
Call this first point $P_{n_1}  = 
(n_1, \log M_{n_1})$. The interval $[P_0, P_1]$
 will form the first side of Newton's polygon.
Let's then turn the half-line
 in the same sense around the point $P_{n_1}  = 
 (n_1, \log M_{n_1})$ until it meets a point 
 $P_{n_2}  = (n_2, \log M_{n_2})$. The interval 
 $[P_1, P_2]$ will form the second side of 
 Newton's polygon, and so on. We thus obtain the
 Newton's polygon of the sequence $\log M_0,\,
 \log M_1, \, \log M_2,\ldots, \log M_n,\ldots $.\\
 $\forall n\in\N $, let us denote by $(\log M_n)^{'}$
 the ordinate corresponding to the 
abscissa $n$ of the
Newton's polygon of the sequence 
$(\log M_n)_{n\in\N}$. The sequence 
$\{(\log M_n)'\}_{n\in\N}$ is the largest convex
sequence whose terms are less than the terms of the 
sequence  $(\log M_n)_{n\in\N}$. We have:
$$ (\log M_n)' =\inf\limits_{k\geq 0,\, 0\leq l\leq n}
(\frac{k\log M_{n-l} + l\log M_{n+k}}{k+l} ).$$
We put $M^c_n =\exp(\log M_n)'$ and 
$ M^c = (M_n^c)_{n\in\N}$, obviousely
$$ M_n^c \leq M_n,\,\,\,\forall n\in\N,$$
and 
$$ M^c_n = \inf\limits_{k\geq 0,\, 0\leq l\leq n}
(M^{\frac{k}{k+l}}_{n-l}, M^{\frac{l}{k+l}}_{n+k} ).$$
It is clear that $C_{M^c}([a, b]) \subset 
C_{M}([a, b])$. 
We give some useful properties of log-convex 
sequences.
\begin{proposition}
    Asuume that the sequence $M=(M_n)_{n\in\in N}$,
   $ M_0=1$, is a log-convex sequence of positive
   real numbers, then 
\begin{enumerate}
    \item [a)] the sequence 
    $(\frac{M_n}{M_{n+1}})_{n\in\N}$ is monotone 
    increasing,
    \item [b)] the sequence  $ (\sqrt[n]{M_n})_{n\in\N}$ 
    is monotone increasing.
\end{enumerate}   
\end{proposition}
\begin{proof}
    a) Is clear from the definition.\\
    b) $$M_n =\frac{M_n}{M_0} = \prod\limits_{j=1}^n
    \frac{M_j}{M_{j-1}} \leq 
    (\frac{M_n}{M_{n-1}})^n$$
    which gives $M_{n-1}^n \leq M_{n}^{n-1}$, or 
    equivalently $\sqrt[n-1]{M_{n-1}} \leq 
    \sqrt[n]{M_n}$.
\end{proof}
Using the convex regularization by means 
of logarithm, we give other conditions equivalent to
Carleman's one, see \cite{Man}.
\begin{theorem}
( Mandelbrojt). Let $M = (M_n )_{n\in\N}$ 
be a sequence of positive real numbers. Then,
    the following conditions are equivalent 
\begin{enumerate}   
  \item $\sum\limits_{n=0}^\infty 
  \frac{1}{\beta_n} = \infty $, where 
  $\beta_n = \inf\limits_{k\geq n}\sqrt[k]{M_k}$
  \item $\sum\limits_{n=0}^\infty 
  \frac{1}{\sqrt[n]{M^c_n}} = \infty $
  \item   $\sum\limits_{n=0}^\infty 
  \frac{M^c_{n-1}}{M^c_n} = \infty $
\end{enumerate}      
\end{theorem}
\begin{proof}
Since the sequence $ (\sqrt[n]{M^c_n})_{n\in\N}$
is increasing, we have
$$ \sqrt[n]{M^c_n} = \inf\limits_{k\geq n}
\sqrt[k]{M^c_k}\leq \inf\limits_{k\geq n}
\sqrt[k]{M_k} =\beta_n,$$
hence 
$$ \sum\limits_{n=0}^\infty
\frac{1}{\sqrt[n]{M^c_n}} \geq 
\sum\limits_{n=0}^\infty
\frac{1}{\beta_n}.$$
Consider the sequence $(N_n )_{n\in\N}$ defined by
$$ N_i = M_{n_i}\,\,\,\mbox{and}\,\,\, N_p = \infty 
\,\,\,\mbox{if}\,\,\, p\neq n_i,\forall i\in\N.$$
where $n_1, n_2 , \ldots  , n_i , \ldots $ 
is the principal sequence of 
$M = (M_n )_{n\in\N}$. Put
$$\gamma_n = \inf\limits_{k\geq n}
\sqrt[k]{N_k} .$$
We see that $\gamma_n \geq \beta_n$, 
$\forall n\in\N$ and if $n\in ]n_{i-1}, n_i]$, we 
have $ \gamma_n = \sqrt[n_i]{M_{n_i}}$.
We then deduce $$ \sum\limits_{n=n_{i-1} +1}^{n_i}
\frac{1}{\gamma_n} = \frac{n_i - n_{i-1}}
{\sqrt[n_i]{M_{n_i}}}.$$
Taking into account that the sequence 
$(\frac{M^c_{n+1}}{M_n^c})_{n\in\N}$  is monotone 
increasing, and for all $n_i$, 
$M_{n_i}^c = M_{n_i}$, we have:
$$ \left(M_{n_i}\right)^{\frac{1}{n_i}} \leq 
\left(\frac{M_{n_i}}{M_{n_{i-1}}}\right)^{\frac{1}
{n_i - n_{i-1}}},$$
hence 
$$ \frac{\log M_{n_i}}{n_i}\leq 
\frac{\log M_{n_i} - \log M_{n_{i-1}}}
{n_i - n_{i-1}}.$$
But we know that the restriction of the function 
$n\mapsto \log M_n^c$  to the interval 
$[n_{i-1}, n_i ]$ is 
linear and
$$ \forall n\in [n_{i-1}, n_i ], \,\, 
\log M_n^c =
\frac{\log M_{n_i} - \log M_{n_{i-1}}}
{n_i - n_{i-1}} n + K $$
for some positive constant $K$. We then obtain
$$ \log M_{n_i} \leq 
\frac{\log M_{n_i} - \log M_{n_{i-1}}}
{n_i - n_{i-1}} 
= \log M_n^c - \log M_{n-1}^c,\,\,\,\forall 
n\in ]n_{i-1}, n_i ].$$ Hence we have:
$$ \sqrt[n_i]{M_{n_i}} \leq \frac{M_n^c}{M_{n-1}^c},
\,\,\,\, \forall 
n\in ]n_{i-1}, n_i ],$$
which gives
$$ \sum\limits_{n=1}^\infty \frac{1}{\gamma_n}
\geq \sum\limits_{n=1}^\infty
\frac{M_{n-1}^c }{M_{n}^c}.$$
Since $\gamma_n \geq \beta_n$, $\forall n\in \N$,
we have 
$$ \sum\limits_{n=1}^\infty \frac{1}{\beta_n}
\geq \sum\limits_{n=1}^\infty
\frac{M_{n-1}^c }{M_{n}^c}.$$
For the last step of the proof, we use Carleman's
 inequality, see \cite[page 112]{Carleman} and 
 \cite{Mal} for the 
 continuous version, 
 $$ a_1 + \sqrt[2]{a_1.a_2} +\sqrt[3]{a_1.a_2.a_3}
 +\ldots + \sqrt[n]{a_1.a_2.a_3 +\ldots . a_n } 
 \ldots \leq e (a_1 + a_2 +\ldots +a_n +  \ldots ),$$
 where $a_1 , a_2 ,\ldots ,a_n,  \ldots$ are real 
 positive  numbers. If we put 
 $ a_n = \frac{M_{n-1}^c }{M_{n}^c}$, we obtain 
 $$ \sum\limits_{n=0}^\infty \frac{1}{\sqrt[n]{M_n^c}}
\leq e \sum\limits_{n=1}^\infty
\frac{M_{n-1}^c }{M_{n}^c}.$$ 
\end{proof}
We need some tools to prove Carleman's theorem.
\section {Some machinery from the theory of 
metric space.}
Let $\mathcal{S}(\R)$  denote the set of all 
real sequences and $\mathcal{P}\subset\N$ an 
infinite set. We will construct a
metric on $\mathcal{S}(\R)$, see \cite{Bang}.
\begin{definition}
    For any $X =(x_n)_{n\in\N}\in \mathcal{S}(\R)$,
    we define 
    $$ \| X\| = \inf\limits_{k\in\mathcal{P} }
    \left ( \max\left ( e^{-k}, \max
    \limits_{0\leq n\leq k}\mid x_n\mid \right)
    \right)$$
\end{definition}
We see that $\| X\|\geq 0$ and $\| X\| = \|- X\|$.\\
If $X =(x_n)_{n\in\N}\in \mathcal{S}(\R)-\{0\}$, 
we can assume that $ x_0\neq 0$  
(if not we change the numbering), we put
$$ \Delta =\{ k\in \mathcal{P}\,/\, e^{-k} < 
\mid x_0\mid \},$$ note that $\Delta\neq \emptyset$,
hence there exits $k'\in \mathcal{P}$ such that 
$\inf \Delta = k'$.
\begin{lemma}
  If $X =(x_n)_{n\in\N}\in \mathcal{S}(\R)-\{0\}$, 
  then there exists  $k\in \mathcal{P}$, 
  $0\leq k\leq k'$,
such that: 
$$ \| X\| = 
 \max\left ( e^{-k}, \max
\limits_{0\leq n\leq k}\mid x_n\mid \right)
$$
\end{lemma}
\begin{proof}
    We can assume that $ x_0\neq 0$.
    We observe that the sequence 
    $\left(\max\limits_{0\leq n\leq k} \mid x_n\mid
    \right)_{k\in\N}$ is increasing in 
    $k$ with infimum $\mid x_0\mid$. As the sequance
    $(e^{-n})_{n\in\N}$ is decresing,
    it follows that for all $k\in \mathcal{P}$, $ k
\geq k'$, we have:
$$ \max\left ( e^{-k}, \max
\limits_{0\leq n\leq k}\mid x_n\mid \right) = 
\max
\limits_{0\leq n\leq k}\mid x_n\mid,$$  
consequently,
$$ \| X\| = \inf\limits_{k\in \mathcal{P},\, 
k\leq k'}\max\left ( e^{-k}, \max
\limits_{0\leq n\leq k}\mid x_n\mid \right).
$$
The result now follows by observing that the infimum
 has been taken over a finite set.
\end{proof}
\begin{lemma}
    If $X =(x_n)_{n\in\N}\in \mathcal{S}(\R)$
    satisfies the inequality:
    $$ e^{-k_1}   \leq \| X\|\leq e^{-k_2},$$
with $k_1,\, k_2\in \mathcal{P}$. Then there exists 
$k\in\mathcal{P}$, $ k_2\leq k\leq k_1$, such that
$$ \| X\| = \max\left ( e^{-k}, \max
\limits_{0\leq n\leq k}\mid x_n\mid \right).$$
\end{lemma}
\begin{proof}
  by lemma 3.2,   there exists  $k_0\in \mathcal{P}$, 
  $0\leq k_0\leq k'$, such that: 
  $$ \| X\| = 
   \max\left ( e^{-k_0}, \max
  \limits_{0\leq n\leq k_0}\mid x_n\mid \right).
  $$
Since $\| X\|\leq e^{-k_2}$, we have $ e^{-k_0}\leq
e^{-k_2}$, hence $ k_2\leq k_0$.\\
If $e^{-k_1}   \leq \| X\|$, then for all $ p\in\N,
\,\, p >k_1$,
$$ e^{-p} < e^{-k_1} \leq \| X\|= 
 \inf\limits_{k\in\mathcal{P} }
    \left ( \max\left ( e^{-k}, \max
    \limits_{0\leq n\leq k}\mid x_n\mid \right)
    \right)$$
    we see then, if $ p >k_1$, we have 
    $$ \max\left ( e^{-p}, \max
    \limits_{0\leq n\leq p}\mid x_n\mid \right)
     = \max
    \limits_{0\leq n\leq p}\mid x_n\mid .$$ 
    It follows that:
    $$ \| X\| = \inf\limits_{k\in\mathcal{P},\,k\leq k_1 }
    \left ( \max\left ( e^{-k}, \max
    \limits_{0\leq n\leq k}\mid x_n\mid \right)
    \right)$$
    which ends the proof.
\end{proof}
\begin{lemma}
    If $X,Y\in \mathcal{S}(\R)$, then we have:
 \begin{enumerate}
    \item $\| X\| \Leftrightarrow  X =0 \,\,\,
    \mbox{(the zero sequence)}$.
    \item $\| X + Y \| \leq \| X\| + \| Y\|$
 \end{enumerate}   
\end{lemma}
\begin{proof}
\begin{enumerate}
 \item  
It is clear that if $ X = 0 $, then 
$\| X  \| =  \inf\limits_{k\in\mathcal{P} }
 e^{-k} = 0 $. \\
To prove the converse, suppose that $ X\neq 0$, 
then there exists $ p\in\mathcal{P}$  such that
$\mid x_p\mid >0$ and for all $q\in \mathcal{P}$, 
$q < p$, we have $x_q =0$. This implies that
$$ \forall j\in \mathcal{P},\,\, j\geq p,\,\, 
\max\left( e^{-j}, \max
    \limits_{0\leq n\leq j}\mid x_n\mid \right)
\geq \mid x_p\mid,$$
and 
$$ \forall j\in \mathcal{P},\,\, j <  p,\,\, 
\max\left( e^{-j}, \max
    \limits_{0\leq n\leq j}\mid x_n\mid \right)
=  e^{-j},
$$
hence $\| X \| \geq \min\left( \mid x_p\mid, 
e^{-p+1}\right) > 0$. 
\item 
If $X =0$ or $Y=0$, it is clear that 
$\| X + Y \| \leq \| X\| + \| Y\|$.\\
Suppose that $\| X\|\neq 0$ and $\| Y\|\neq 0$. By 
lemma 3.2, we have 
$$\| X\| = 
 \max\left ( e^{-k_1}, \max
\limits_{0\leq n\leq k_1}\mid x_n\mid \right),$$
$$\| Y\| = 
 \max\left ( e^{-k_2}, \max
\limits_{0\leq n\leq k_2}\mid x_n\mid \right),$$
for some $k_1,\ k_2 \in \mathcal{P}$. We can 
suppose that $ k_1\leq k_2$. For all $ n\leq k_1$,
we have 
$$ \mid x_n + y_n\mid \leq \mid x_n\mid + 
\mid x_n\mid \leq \max
\limits_{0\leq n\leq k_1}\mid x_n\mid +
\max
\limits_{0\leq n\leq k_2}\mid y_n\mid \leq \| X\|+
\| Y\| ,$$
hence $$ \max
\limits_{0\leq n\leq k_1}\mid x_n + y_n\mid \leq 
\| X\|+\| Y\|.$$
Note also that $$e^{k_1}\leq \| X\|+\| Y\|, $$
hence 
$$ \| X + Y \|\leq \max\left ( e^{-k_1}, \max
\limits_{0\leq n\leq k_1}\mid x_n + y_n \mid \right)
\leq \| X\|+\| Y\|.$$
\end{enumerate}
\end{proof}
We can then provide the space $\mathcal{S}(\R)$ 
with a distance function, $d$, 
defined by $ d(X,Y) = \| X - Y \|$.\\

Having developed this machinery, we now utilize 
it in the proof of the statement 
that condition
iii) of theorem 2.7 implies quasi-analyticity 
of the relevant functions.\\
Let $M =(M_n)_{n\in N}$ be a sequence of 
positive numbers with $M_0 =1$, and suppose 
that $\liminf\limits_{n\to \infty}
\sqrt[n]{M_n(f)} = \infty$. Consider 
$M^c = (M_n^ c)_{n\in \N}$ the convex 
regularization by means of the logarithm of the 
sequence $M =(M_n)_{n\in N}$. We denote by 
$\mathcal{P}\subset \N$ the set of all $p\in\N$,
such that $M_p = M_p^c$, we suppose that 
$0\in \mathcal{P}$.\\
If $f\in C_M([a,b])$,  $ \forall t\in [a,b]$, we put
$$ X_f(t)= (x_{f,n}(t))_{n\in \N},\,\,\mbox{where}\,\,
x_{f,n}(t) = \frac{f^{(n)}(t)}{M_n^c e^n},\,\,
\forall n\in\N$$
and 
\begin{equation}
\| X_f(t)\| = \inf\limits_{k\in\mathcal{P} }
    \left ( \max\left ( e^{-k}, \max
    \limits_{0\leq n\leq k}\mid x_{f,n}(t)\mid \right)
    \right).
\end{equation}    
    By lemma 3.2, there exists $l\leq k'$, 
    $ l\in \mathcal{P}$, such that
    $$ \| X_f(t)\| = \max\left ( e^{-l}, \max
    \limits_{0\leq n\leq l}\mid x_{f,n}(t)\mid 
    \right),$$
    where $k'$ is the smallest $p\in \mathcal{P}$,
    such that $e^{-p} < \mid f(t)\mid$. We remark that:
    $$ e^{-l}\leq \| X_f(t)\|, \,\,\,\mbox{and }
    \,\,\,\frac{\mid f^{(n)}(t)\mid}{M_n^c e^n}\leq \| X_f(t)\|,
    \,\,\forall n= 0, 1,\ldots, l.$$
    Let $\tau\in\R$ such that $ t+\tau\in [a,b]$.
    By (3.1), we have 
    \begin{equation}
      \| X_f(t+\tau )\| \leq  
    \max\left ( e^{-l}, \max
    \limits_{0\leq n\leq l}
    \frac{\mid f^{(n)}(t+\tau)\mid}{M_n^c e^n} 
    \right).
    \end{equation}
    The following lemma gives us a link between
    $\| X_f(t+\tau )\|$ and $\| X_f(t)\|$.
\begin{lemma}
    Suppose that $\| X_f(t)\|\neq 0$, then 
    $$ \| X_f(t+\tau )\|\leq \| X_f(t)\| \exp\left(
        e\mid \tau\mid \frac{M^c_l}{M^c_{l-1}}
    \right)
    $$
    where $l$ satisfies
    $$\| X_f(t)\| = \max\left ( e^{-l}, \max
    \limits_{0\leq n\leq l}
    \frac{\mid f^{(n)}(t)\mid }{M_n^c e^n} 
    \right).$$
    
\end{lemma}
\begin{proof}
According to above, there exists $l\in 
\{ 1,2,\ldots,k'\}\cap \mathcal{P}$ such that
$$\| X_f(t)\| = \max\left ( e^{-l}, \max
\limits_{0\leq n\leq l}
\frac{\mid f^{(n)}(t)\mid }{M_n^c e^n} 
\right).$$
where $k'$ is the smallest $j\in  \mathcal{P}$ 
such that $ e^{-j} < \mid f(t)\mid $.\\
If $e^{-l}\geq \max
\limits_{0\leq n\leq l}
\frac{\mid f^{(n)}(t+\tau)\mid }{M_n^c e^n}$, 
we have, by (3.2),
$$ \| X_f(t+\tau)\| \leq e^{-l} \leq \| X_f(t)\|,$$
and hence the statement of the lemma holds true.\\
Suppose not, then there exists $0\leq n\leq l$,
such that $ \| X_f(t+\tau)\|\leq 
\frac{\mid f^{(n)}(t+\tau)\mid }{M_n^c e^n}$.
By using Taylor's
formula  for the function $f^{(n)}$ at the point 
$t$, we get
\begin{align*}
\mid f^{(n)}(t+\tau)\mid  &\leq  
\sum\limits_{j=0}^{l-n-1}\frac{\mid\tau\mid^j}{j!}
\mid f^{(n+j)}(t)\mid + 
\frac{\mid\tau\mid^{l-n}}{(l-n)!}
f^{(l)}(\xi) \\
&\leq  \sum\limits_{j=0}^{l-n-1}
\frac{\mid\tau\mid^j}{j!}M_{n+j}^c e^{n+j}\| X_f(t)\|
+ \frac{\mid\tau\mid^{l-n}}{(l-n)!}
M_{l}^c e^{l}\| X_f(t)\|\\
        &= \| X_f(t)\|\sum\limits_{j=0}^{l-n}
        M_{n+j}^c e^{n+j}.
\end{align*}
Then 
\begin{align*}
\frac{\mid f^{(n)}(t+\tau)\mid}{M_n^ce^n} &\leq 
\| X_f(t)\|\sum\limits_{j=0}^{l-n}
\frac{\mid\tau\mid^j}{j!}
\left(\frac{M_{n+j}^c}{M_n^c}
\right)e^j\\&\leq 
\| X_f(t)\|\sum\limits_{j=0}^{l-n}
\frac{\mid\tau\mid^j}{j!}
\left(\frac{M_{l}^c}{M_{l-1}^c}
\right)^je^j\\ & \leq \| X_f(t)\|
\exp\left( e\mid \tau\mid \frac{M_{l}^c}{M_{l-1}^c}
\right)
\end{align*}
where log-convexity of the sequence 
$(M_n^c)_{n\in\N}$ has 
been used to derive $ \frac{M_{n+j}^c}{M_n^c}\leq 
\left(\frac{M_{l}^c}{M_{l-1}^c}\right)^j$. Hence the
lemma is proved. 
\end{proof}
\section{Monotonicity property for quasianalytic 
Denjoy-Carleman class}
We know that an analytic function $f$ on the 
interval $[a, b]$ is entirely determined by 
the element
$\left({f^{(n)} (c)}\right)_{n\in\N}$ , where $c \in [a, b]$. 
We are interested in a generalization of this fact,
 which can be stated as follows:\\
 Let $(x_n )_{n\in\N}$ be a sequence of elements of 
 the interval $[a, b]$. Which condition must check 
 the sequence $(x_n )_{n\in\N}$ in order that the 
 element $(f^{(n)}(x_n) )_{n\in\N}$ determines the 
 function $f $ completely. We see that
 if the sequence $(x_n )_{n\in\N}$ is constant 
 we get our first property.
 For an analytic function, a response is given 
 by W. Gontcharoff \cite{Gon}.
\begin{theorem}
Let $f$  be an analytic function on the interval 
$[a, b]$. The function $f $ is entirely
determined by the knowledge of the values 
$(f^{(n)}(x_n) )_{n\in\N}$ if the series
$\sum\limits_{n=1}^\infty \mid x_{n-1} - x_n\mid$
converges.   
\end{theorem}
As a consequence, we deduce that if 
$f^{(n)}(x_n) =0$, $\forall n\in\N$, and the series
$\sum\limits_{n=1}^\infty \mid x_{n-1} - x_n\mid$
converges, then the function f is identically zero.

The question now is whether a similar result remains
 valid for a quasianalytic class. \\ In the case
of a quasianalytic class of Denjoy-Carleman, 
we have the following theorem proved by W.Bang
[1].\\
Let $M = (M_n )_{n\in\N }$ be a sequence of real
positive numbers with $M_0 = 1$. Suppose 
that $M = (M_n )_{n\in\N }$ is
logarithmically convex and satisfying 
one of the equivalent conditions of Proposition 2.6.
In other words, the class is quasianalytic.
\begin{theorem}
 Let $f\in C^\infty ([a,b])$ such that 
 $$ \sup\limits_{t\in[a,b]}\mid f^{(n)}(t)\mid 
 \leq M_n,\,\,\,\,\forall n\in\N.$$
 Suppose that there exists a sequence 
 $(x_n )_{n\in\N}$, 
  $x_n\in [a, b],\,\forall n\in\N$, 
  such that $f^{(n)}(x_n ) = 0, \,\,\,
  \forall n\in\N$. If f is non-identically zero, 
  then the series 
  $\sum\limits_{n=1}^\infty \mid x_{n-1} - x_n\mid$
  is divergent.
  \end{theorem}
For the convenience of the reader and for completeness, 
we reproduce the proof.\\
First, we introduce some preliminaries and 
definitions.\\ We consider a  not identically zero
function $f$ which 
satisfies the hypotheses of the Theorem 4.2.
For each $n\in\N$, 
we define a  function 
 $[a,b]\ni t\mapsto B_{f,n}(n)$,  by 
 $$ B_{f,n}(t) = \sup\limits_{j\geq n}
 \frac{\mid f^{(j)}(t)\mid } {e^j M_j}
.$$
 First, we list some properties of this sequence of 
 functions. 
 \begin{enumerate}
    \item $\forall t\in[a,b], \,\,\,\,\,\,
    B_{f,n}(t) \leq e^{-n}$,
    \item $\forall t\in[a,b], \,\,\,\,\,\, 
    B_{f,0}(t)\geq B_{f,1}(t)\geq B_{f,2}(t)
    \geq \ldots B_{f,n}(t)\geq \ldots $
    \item if $f^{(n)}(t_0)=0$, then $B_{f,n}(t_0) =
    B_{f,n+ 1}(t_0)$.
  \item Since f is non-identically zero, we have, 
  $\forall t\in[a,b]$, $B_{f,n}(t)\neq 0$.\\
 Indeed, if there exists $t_0\in [a,b]$ such that
 $B_{f,n}(t_0) = 0$. As $f^{(n)}$ is in the same 
 quasianalytic class as $f$, we deduce that $f^{(n)}$ is 
 identically zero. Since $ f^{(n-1)}(x_{n-1}) =0$, 
 the function $ f^{(n-1)}$ is also identically zero.
 Continuing, we find that the function $f$ is 
 identically zero, which is a contraduction.
\end{enumerate}
 \begin{lemma}
    The function $ [a,b]\ni t\mapsto B_{f,n}(t)$ 
    satisfies the estimate
    $$ B_{f,n}(t+\tau)\leq \max 
\left(B_{f,n}(t), e^{-q}\right) 
\exp\left(e\mid \tau\mid \frac{M_q}{M_{q-1}}\right),$$
$\forall q\in\N$, $ q > n$ and forall $ t,\tau +t\in 
[a,b]$.
\end{lemma}
\begin{proof}
    We follow the proof of the lemma 3.5. 
    Let $j\in\N$, $n\leq j < q$,
\begin{align*}   
\frac{f^{(j)}(t+\tau)}{e^j M_j}&\leq \sum\limits_{i=0}
^{q-j-1}\frac{\mid \tau\mid^i}{i!e^j M_j}
\mid f^{(i+j)}(t)\mid  + 
\frac{\mid f^{(q)}(\xi)\mid}{(q-j)!e^j M_j}
\\ & = \sum\limits_{i=0}^{q-j-1}\frac{M_{i+j}}{M_j}
\frac{f^{(i+j)}(t)}{e^{i+j} M_{i+j}}
\frac{(e\mid\tau\mid)^i}{i!}+e^{-q}\frac{M_q}{M_j}
\frac{\mid f^{(q)}(\xi)\mid }{M_q}
\frac{(e\mid\tau\mid )^{q-j}}{(q-j)!}\\
&\leq B_{f,n}(t)\sum\limits_{i=0}^{q-j-1}
\frac{(e\mid\tau\mid)^i}{i!}
\left(\frac{M_q}{M_{q-1}}\right)^i + e^{-q}
\left(\frac{M_q}{M_{q-1}}\right)^{q-j}
\frac{(e\mid\tau\mid )^{q-j}}{(q-j)!}\\
&\leq \max \left( B_{f,n}(t), e^{-q}\right)
\exp\left( e\mid \tau\mid \frac{M_q}{M_{q-1}}\right).
\end{align*} 
Where we used that the sequence $ M= (M_n)_{n\in\N}$
is logarithmically convex and   
$\frac{\mid f^{(q)}(\xi)\mid }{M_q}\leq 1$.
\end{proof}
As a consequence, we see that the function 
$ [a,b]\ni t\mapsto B_{f,n}(t)$  is continuous. 
Indeed, if $ t_0\in[a,b]$, since $B_{f,n}(t_0) >0$, 
there exists $q\in\N$, $ q >n$, such that 
$ e^{-q} \leq B_{f,n}(t_0)$. By lemma 4.3., we have 
$$ \mid B_{f,n}(t_0+\tau) - B_{f,n}(t_0)\mid \leq 
B_{f,n}(t_0)
\left( \exp \left(e\mid \tau \mid 
\frac{M_q}{M_{q-1}}\right) -
1\right).$$
Let $f$ and $(x_n)_{n\in\N}$ as in the theorem.\\
For $k \geq 1$, 
set $\tau_k =\sum\limits_{n=0}^{k-1}
\mid x_n - x_{n+1}\mid $,  and 
$\tau_0= 0$. If $t\in [\tau_{n-1},\tau_n]$, then

$$\left\lbrace
\begin{array}{ll}
x_{n-1}+\tau_{n-1}-t \in [x_n, x_{n-1}]\subset 
[a,b] & 
\mbox{if } x_n < x_{n-1}\\
x_{n-1}-\tau_{n-1} +t \in [x_{n-1}, x_n]\subset 
[a,b] & \mbox{if } x_{n-1} < x_n
\end{array}
\right.$$
We define a function $\overline{B}_{f,n}$ on 
$[\tau_{n-1},\tau_n]$ by
$$\overline{B}_{f,n}(t)=\left\lbrace
\begin{array}{ll}
B_{f,n}(x_{n-1}+\tau_{n-1} -t) & 
\mbox{if }  x_n < x_{n-1}\\
B_{f,n}(x_{n-1}-\tau_{n-1} +t) & \mbox{if } 
 x_{n-1} < x_n
\end{array}
\right.$$
The function $t\mapsto \overline{B}_{f,n}(t)$ is
continuous and the property 3)
 $$ \overline{B}_{f,n}(\tau_n) = B_{f,n}(x_{n}) =
 B_{f,n+1}(x_{n}) = \overline{B}_{f,n+1}(\tau_n).$$
 We can then paste together the functions $
 [\tau_{n-1},\tau_n]\ni t\mapsto 
 \overline{B}_{f,n}(t)$ with different $n\in\N$, and 
 define a new function, $t\mapsto b_f(t)$,
  on the interval $[0,\tau[$,
 where $ \tau =\sup\limits_{n\in\N}\tau_n$, by 
 $$ b_f(t)= \overline{B}_{f,n}(t),\,\,\,\,\mbox{if}
 \,\,\,\,t\in [\tau_{n-1},\tau_n].$$
 This is a continuous function. By 1. and 2., 
 we find that
$$ b_f(t) \leq e^{-n},\,\,\,\,
 \forall t \geq \tau_{n-1},\leqno{(*)}$$
and hence $b_f(t) \to 0$ when  $ t\to \tau$.\\
Note that for all $t_1, \,t_2 \in 
[\tau_{n-1},\tau_n]$, we have 
\begin{equation}
    \forall q\geq n,\,b_f(t_2)\leq \max \left(
        b_f(t_1), e^{-q}\right)\exp \left( e\mid
        t_1 -t_2\mid \frac{M_q}{M_{q-1}}\right).      
\end{equation}
Indeed, since $t_1, \,t_2 \in 
[\tau_{n-1},\tau_n]$, we have 
$$b_f(t_2) = \overline{B}_{f,n}(t_2)=
\left\lbrace
\begin{array}{ll}
B_{f,n}(x_{n-1}+\tau_{n-1} -t_2) & 
\mbox{if }  x_n < x_{n-1}\\
B_{f,n}(x_{n-1}-\tau_{n-1} +t_2) & \mbox{if } 
 x_{n-1} < x_n
\end{array}
\right.$$
Using Lemma 4.3, we deduce (4.1).\\
{\it{Proof of Theorem 4.2.}}\\
Since $f$ is non-identically zero, the function
$t\mapsto b_f(t)$ is also non-identically zero.
Then, there exists $\mu >0$, such that 
$]0,\mu] \subset Im \,(b_f)$. Let $k_0\in\N $ be the 
smallest integer such $ e^{-k_0}\in ]0,\mu] 
\subset Im \,(b_f)$. We put
$$ t_{k_0}= \inf\{ t\in [0,\tau[\,/\, b_f(t) 
= e^{-k_0}\},$$
and for $ k> k_0$
$$ t_{k}= \inf\{ t\in ]t_{k-1},\tau[\,/\, b_f(t) 
= e^{-k}\}.$$
So, we have a strictly increasing sequence 
$(t_k)_{k\geq k_0}$ such that $b_f(t_k) = e^{-k}$.
We see that, for each $ t\in]t_{k-1}, t_k[$, 
$b_f(t) > e^{-k}$  
 and by $(*)$,
$ t_k < \tau_k,\,\,\forall k \geq k_0$.
The sequence $(\tau_n)_{n\in\N}$ defines on each
interval $]t_{k-1}, t_k[$, $k\geq k_0$,  
a subdivision. \\
Let $s\in\N$  be the largest integer such that 
$ \tau_s \leq t_{k-1}$ and $r\in\N$ 
the smallest integer such that 
$ t_{k} \leq \tau_r$.  Using
the fact that for each $ t\in ]t_{k-1}, t_k[$, we 
have $b_f(t) > e^{-k}$ and (4.1), we obtain

$$
\left\{
    \begin{array}{ll}
        b_f(t_{k-1}) \leq b_f(\tau_{s+1})  
        \exp \left( e
        (\tau_{s+1} -t_{k-1}) \frac{M_k}{M_{k-1}}
        \right)
         \\
         b_f(\tau_{s+1}) \leq b_f(\tau_{s+2})  
        \exp \left( e
        (\tau_{s+2} -\tau_{s+1}) \frac{M_k}{M_{k-1}}
        \right) \\
        \ldots \ldots  \ldots \ldots \ldots \ldots
        \ldots \ldots \ldots \ldots \ldots \ldots
        \ldots \ldots \\
        \ldots \ldots  \ldots \ldots \ldots \ldots
        \ldots \ldots \ldots \ldots \ldots \ldots
        \ldots \ldots \\
        b_f(\tau_{r-1}) \leq b_f(t_{k}) \,\, 
        \exp \left( e
        (t_{k} -\tau_{k-1}) \frac{M_k}{M_{k-1}}
        \right) 
    \end{array}
\right.
$$
We infer that
$$ b_f(t_{k-1}) \leq b_f(t_{k}) \,\, 
\exp \left( e
(t_{k} -t_{k-1}) \frac{M_k}{M_{k-1}}
\right).
$$
Since $b_f(t_k)= e^{-k}$, $\forall k\geq k_0$, we
have $$ t_{k} -t_{k-1} \geq \frac{1}{e}
\frac{M_{k-1}}{M_k},\,\,\,\,\forall k\geq k_0,$$
hence 
$$ t_k \geq t_{k_0} + \frac{1}{e}
\sum\limits_{j=k_0+1}^k\frac{M_{j-1}}{M_j}.$$
Since $ t_k < \tau_k,\,\,\forall k \geq k_0$, we 
 obtain
$$ \sum\limits_{j=0}^{k-1}\mid x_j -x_{j+1}\mid  >
\frac{1}{e}
\sum\limits_{j=k_0+1}^k\frac{M_{j-1}}{M_j},$$
which proves the theorem.\\
As consequence of  Theorem 4.2, we have the 
following result:\\
Let $M = (M_n )_{n\in\N}$ be a sequence of 
positive real numbers with $M_0 = 1$. 
Suppose that $M = (M_n )_{n\in\N}$ is
logarithmically convex and satisfying one 
of the equivalent conditions of Theorem 2.7.
\begin{corollary}
    Let $ f\in C^\infty([a,b])$, such that 
    $$\sup\limits_{t\in [a,b]}\mid f^{n}(t)\mid \leq
    M_n, \,\,\,\,\forall n\in \N.$$
    If $f^{n}(a) > 0$, $\forall n\in \N$, then 
    $f^{n}(x) > 0$, $\forall n\in \N,\,\,\, \forall
    x\in [a,b]$.
\end{corollary}
\begin{proof}
    Suppose there exits $k_0\in\N$, such that 
    $f^{(k_0)}$ has a zero $x_{k_0}\in ]a,b]$. Then
    there exists $ x_{k_0+1} < x_{k_0}$, such that  
    $f^{(k_0+1)}(x_{k_0+1}) = 0$.
    Continuing, we find a strictly decreasing 
    sequence 
$x_{k_0} > x_{k_0+1} > . . .$  where $x_l$
 is a zero of $f^{(l)} $, 
for all $ l\geq  k_0$. By
Theorem 4.2, the function $f ^{(k_0)}$ is 
identically $0$, which contradicts $f ^{(k_0)}(a) >
0$.
\end{proof}
As an immediate consequence of  Corollary 4.4., 
we deduce, according to Bernstein's theorem,
see \cite[page 146]{Wider},  that the function $f$
 can be extended analytically into the 
 plane of complex
numbers to a holomorphic function in 
the disk 
$ \{z\in\C\,/\,  \mid z- a\mid  < b - a\}$.
\section{ Borel mapping }
The Borel mapping takes germs, at the origin in 
$\R$, of $C^\infty$ functions, $\mathcal{E}$, 
 to the sequence of iterated derivatives at 
 the origin. It is a classical result due to 
 Borel that the Borel mapping is surjective and 
 not injective.
A subring $\mathcal{B}\subset \mathcal{E}$
 is called quasianalytic if the restriction of the 
 Borel mapping to $\mathcal{B}$ is
injective. As a direct application of the 
corollary 4.4., we give a proof of the  
Carleman result concerning the non-surjectivity of 
the Borel mapping restricted to the subring
of smooth functions in a quasianalytic 
Denjoy-Carleman class which strictly contains the
analytic class.
\begin{theorem} (Carleman). 
Let $C_M([a, b])$ be a quasianalytic class which 
contains strictly the analytic class. 
Then the Borel mapping
\[
T_c : C_M([a, b]) \to R[[x]], \quad c \in [a, b], 
\quad T_c(f) = \sum_{n=0}^{\infty} \frac{f^{(n)}(c)}{n!} x^n.
\]
is not surjective.
\end{theorem}
\begin{proof}
     Consider a formal series with radius of 
     convergence equal to zero,
     $\sum\limits_{n\in\mathbb{N}} a_n x^n$, and 
     suppose that  
      $a_n > 0$, $\forall n \in \mathbb{N}$. 
     By corollary 4.4., we see 
     that $\sum_{n\in\mathbb{N}} a_n x^n 
    \notin T_c (C_M([a, b]))$.
\end{proof}
We can then ask ourselves the following question:\\
Is Carleman's theorem true for any 
quasianalytic subring of $\mathcal{E}$?\\
Here $\mathcal{E}$ is the ring of germs, 
at the origin of $\R$, 
of the  $ C^\infty$ functions.\\
In other words, if $\mathcal{A}\subset \mathcal{E}$
 is a general quasianalytic subring, is the Borel
mapping, $T_0 : \mathcal{A} \to R[[x]],\,\,\,
T_0(f) = \sum\limits_{n=0}^{\infty} 
\frac{f^{(n)}(0)}{n!} x^n$, nonsurjective? \\
The answer is negative, there exist 
quasianalytical subrings for which the Borel 
mapping  is surjective. In other words, Theorem 5.1. 
is not true for general quasianalytic rings, see 
\cite{Elk}.
We consider the following set, ordered by inclusion:
$$\mathcal{F}=\{ \mathcal{A}\subset \mathcal{E}\,/\,
\mathcal{A}\,\,\mbox{ is a quasianalytic subring
closed under differetiation}\,\,\}$$
The following theorem is proved in \cite{Elk}
\begin{theorem}
    Let $\mathcal{A}\subset \mathcal{E}$ be 
    a quasianalytic subring closed under 
    differentation. Then
    $$ \mathcal{A} \,\,\mbox{ is a maximal element
    of}\,\, \mathcal{F} \Longleftrightarrow \,\,
    \mbox{the Borel map}\,\, \mathcal{A}\to
    \R[[x]]\,\,
    \mbox{is surjective}
    $$
\end{theorem}
This theorem shows that  Carleman's 
theorem on the nonsurjectivity of Borel 
mapping for the quasianalytic classes of 
Denjoy–Carleman is not true for certain general 
quasianalytic subrings.\\
Corollary 4.4. provides us with a criterion to 
see if  Theorem 5.1. is true for a 
quasianalytic class.
\section{ Quasi-analytic classes associated to a 
sequence of integers}
For a function $f$ which is $C^\infty$ on $[a, b]$,
 we can formulate the principle of Pringsheim 
 as follows:
\begin{equation}
f \text{ is analytic on the interval } [a, b] 
\Leftrightarrow \limsup_{n\to\infty} \frac{1}{n} 
\sqrt[n]{M_n(f)} < \infty . 
\end{equation}
To preserve the validity of the condition 
$\mathcal{DC})$ 
for a function $f$, one only needs a weakened 
version of the principle of Pringsheim, 
namely the following condition:
\begin{equation}
\liminf_{n\to\infty} \frac{1}{n} 
\sqrt[n]{M_n(f)} < \infty . 
\end{equation}
More precisely we have:
\begin{proposition}
    Let $f$ be a $C^\infty$ function on $[a, b]$ 
    such that
\[
\liminf_{n\to\infty} \frac{\sqrt[n]{M_n(f)}}{n} < 
\infty.
\]
If $c \in [0, 1]$ is such that $f^{(n)}(c) = 0, 
\forall n \in \mathbb{N}$, then $f$ is identically 
null on $[a, b]$.
\end{proposition}
\begin{proof}
    By Taylor's formula, we have, for each 
    $n \in \mathbb{N}$,
    \[
    f(x) = f^{(n)}(c + \theta_n(x - c)) 
    \frac{(x - c)^n}{n!},
    \]
    where $0 < \theta_n < 1$. Let $(n_k)_k$ be 
    an infinite subsequence such that:
    \[
    \lim_{k\to\infty} 
    \frac{\sqrt[n_k]{M_{n_k}(f)}}{{n_k}} = 
    \liminf_{n\to\infty} \frac{\sqrt[n]{M_n(f)}}{n}.
    \]
    Hence, there exists $A > 0$, such that 
    \[
M_{n_k}(f) \leq A^{ n_k} n_k !, \quad \forall n_k.
\]
We have then
\[
f^{(n_k)}(c + \theta_{n_k}(x - c)) \leq 
A^{ n_k} n_k !, \quad \forall n_k,
\]
and consequently
\[
|f(x)| \leq (A|x - c|)^{n_k}, \quad 
\forall k \in \mathbb{N}.
\]
Hence if $|x - c| < \frac{1}{A}$, it follows 
that $f(x) = 0$.
Choosing in place of $c$ the point 
$c \pm \frac{1}{2A}$ and once more repeating 
the same reasoning, we obtain
$f(x) = 0$ on the whole interval $[a, b]$.
\end{proof}
Let $\overline{n} = (n_k)_{k\in\N}$ be a strictly 
increasing 
sequence of natural numbers. We denote 
by $C_{\overline{n}}([a, b])$ the set of 
all $C^\infty$ 
functions on $[a, b]$ such that, there exist 
two positive constants $A$, $B$, such that
\[
M_{n_k}(f) \leq BA^{n_k} n_k!, 
\quad \forall k \in \mathbb{N}.
\]
By proposition 6.1., any sequence of natural numbers 
$\overline{n} = (n_k)_{k\in\N}$, which increases 
without limit, defines some quasianalytic class 
of functions which satisfies the condition 
$\mathcal{DC})$. 
We will call such class, quasi-analytic class 
with respect to the sequence 
$\overline{n} = (n_k)_{k\in\N}$.
\subsection {Monotonicity property for 
quasianalytic classes associated to a 
sequence of integers}
\subsubsection{The Abel-Gontcharoff polynomials}
In this subsection, we review some known results
about Abel-Gontcharoff polynomials that can be found
in \cite{Gon}.\\
Let $(x_n)_{n\in\N}$ be a sequence of real numbers.
\begin{definition}
     The Abel-Gontcharoff polynomial of degree 
    $0$ is $Q_0(x) = 1  $. The Abel-Gontcharoff
    polynomial of degree $n \geq 1$  is defined as
    \[
        Q_{n}(x, x_0, x_1, \dots, x_{n-1}) = 
    \int_{x_0}^x dt_1 \int_{x_1}^{t_1} dt_2 \dots 
    \int_{x_{n-1}}^{t_{n-1}} dt_n.
    \]
\end{definition} 
Therefore, $Q_n,\,\,n\geq 1$, depends only on 
the first $n$ parameters $x_0, x_1,\ldots,x_{n-1}$
of the sequence $(x_n)_{n\in\N}$.\\
 If $(x_n)_{n\in\N}$
is the constant sequence, $x_n = x_0,\,\,\forall 
n\in\N$, we see that $$ Q_{n}(x, x_0, x_0, \dots, x_{0}) =
\frac{(x-x_0)^n}{n!}.$$
From the definition, we easily see 
that
\begin{equation}
\frac{d Q_n(x, x_0,x_1, \dots, x_{n-1})}{dx} =
Q_{n-1}(x, x_1, \dots, x_{n-1})  
\end{equation} 
 and 
 \begin{equation}
 Q_n(x_0, x_0, x_1, \dots, x_{n-1})=0.
\end{equation}
As an immediate consequence of (6.3), we have 
$$ 
Q^{(k)}_n(x, x_0,x_1 \dots, x_{n-1})=
Q_{n-k}(x, x_k, \dots, x_{n-1}),\,\,\,\,
0\leq k\leq n,$$
and from (6.4), 
\begin{equation}
Q^{(k)}_n(x_k, x_0,x_1 \dots, x_{n-1})= 0,
\,\,\,\,0\leq k\leq n-1.
\end{equation}
The system 
$Q_n(x, x_0, x_1, \dots, x_{n-1})$, $n = 0, 1, 2, 
\dots$ is the Abel-Gontcharoff polynomials
associated to the sequence 
$(x_n)_n$.\\
 Using (6.3) and (6.4), Abel-Gontcharoff 
polynomials  can be 
 calculated step by step:
\begin{align*}
    & Q_1(x, x_0) = (x - x_0)\\
    & Q_2(x, x_0, x_1) = \frac{1}{2} \left((x - x_1)^2 - (x_0 - x_1)^2\right)\\
    & Q_3(x, x_0, x_1, x_2) = \frac{1}{3!} 
    \left((x - x_2)^3 - 3(x_1 - x_2)^2(x - x_0) - 
    (x_0 - x_2)^3\right)
    \end{align*}
Let us remark that, from the definition, we have:
\begin{align*}
&\frac{\partial Q_n(x, x_0, x_1, \dots, x_{n-1}) }
    {\partial x_k}  = \\ & - 
    Q_k(x, x_0, x_1, \dots, x_{k-1})
    Q_{n-k-1}(x_k, x_{k+1}, \dots, x_{n-1})  
\end{align*}   
\begin{proposition}
 Let  $Q_n(x, x_0, x_1, \dots, x_{n-1})$ be the 
 degree $n\geq 1$  Abel-Gontcharoff polynomial, 
 $y$, a new variable and let $0\leq k <n$. Then
 \begin{align*} 
 & Q_n(x, x_0, x_1, \dots, x_{n-1}) -  
 Q_n(x, x_0, x_1, \dots, x_{k-1}, y, x_{k+1}, 
 \ldots, x_{n-1})  \\ = & 
 Q_k(x, x_0, x_1, \dots, x_{k-1})
 Q_{n-k}(y, x_k, \dots, x_{n-1})
\end{align*}
\end{proposition}
\begin{proof}
    The result follows by integrating, 
    \begin{align*} 
        & Q_n(x, x_0, x_1, \dots, x_{n-1}) -  
        Q_n(x, x_0, x_1, \dots, x_{k-1}, y, x_{k+1}, 
        \ldots, x_{n-1}) \\ =& \int_{y}^{x_k}  
    \frac{ 
    \partial Q_n(x, x_0, \dots, x_{k-1}, x_k, x_{k+1}, 
    \ldots, x_{n-1}) }
    {\partial x_k} dx_k \\ =& -\int_{y}^{x_k} 
    Q_k(x, x_0, x_1, \dots, x_{k-1})
    Q_{n-k-1}(t, x_{k+1}, \dots, x_{n-1})dt
\\ =& -Q_k(x, x_0, x_1, \dots, x_{k-1})
\int_{y}^{x_k} Q_{n-k-1}(t, x_{k+1},  
\dots, x_{n-1})dt\\ = &
Q_k(x, x_0, x_1, \dots, x_{k-1})
Q_{n-k}(y,x_k,x_{k+1}, \dots, x_{n-1})
\end{align*}
since by  (6.3), we have 
$$
\frac{d Q_{n-k}(t,x_k,x_{k+1}, \dots, x_{n-1})}
{dt} =
Q_{n-k-1}(t, x_{k+1},  
\dots, x_{n-1}),  
$$      
and   
$$Q_{n-k}(x_k,x_k,x_{k+1}, \dots, x_{n-1}) =0.$$
\end{proof}
\begin{proposition}
    Let  $Q_n(x, x_0, x_1, \dots, x_{n-1})$ be the 
    degree $n\geq 1$  Abel-Gontcharoff polynomial, 
    $y_0, y_1,\ldots,y_{n-1}$, a new variables. Then
\begin{align*}  
    Q_n(x, x_0, x_1, \dots, x_{n-1}) &=
    Q_n(x, y_0, y_1, \dots, y_{n-1}) \\ &+ 
    \sum\limits_{i=0}^{n-1}
    Q_i(x, y_0, y_1, \dots, y_{i-1})
    Q_{n-i}( y_i, x_i, \dots, x_{n-1}) 
\end{align*} 
\end{proposition}
\begin{proof}
Applying Proposition 6.3. for $ 0\leq i <n$, 
we get 
\begin{align*}
Q_n(x, y_0,  \dots, y_{i-1},x_i,\ldots,x_{n-1})-&
Q_n(x,y_0, \dots, y_{i-1},y_i,x_{i+1},
\ldots,x_{n-1})\\ = &
Q_i(x, y_0, y_1, \dots, y_{i-1})
Q_{n-i}( y_i, x_i, \dots, x_{n-1}).
\end{align*}
By taking the sum of the two members of 
the previous equality, we obtain the result.
\end{proof}
Taking $ y_0 =\ldots =  y_{n-1} = 0$ 
in the  Proposition 6.4., we recover the 
standard form of the Abel-Gontcharoff polynomial:
$$ Q_n(x, x_0, x_1, \dots, x_{n-1}) =
\frac{x^n}{n!} + \sum\limits_{i=0}^{n-1}
\frac{x^i}{i!}Q_{n-i}( 0, x_i, \dots, x_{n-1}).
$$
\subsection
{Some estimate of Abel-Gontcharoff polynomial}
An important question is:\\
If $Q_n(x, x_0, x_1, \dots, x_{n-1})$ is the 
$n\geq 1$ degree Abel-Gontcharoff polynomial,
what is the largest value of 
$\mid Q_n(x, x_0, x_1, \dots, x_{n-1})\mid$. In 
\cite{Levinson} it was proved that
$$ \sup\limits_{\mid x\mid \leq 1,\,\mid x_j\mid 
\leq 1,\, o\leq j\leq n-1 }
\mid Q_n(x, x_0, x_1, \dots, x_{n-1})\mid 
\leq r^{n+1}, $$
where $ 0 < r < 1, 386$.\\
This result leads to estimate a value for 
Whittaker's constant, $W$.\\
The Whittaker's constant,$W$, is defined as the 
least upper bound of numbers $c$ such that, if
$f$ is an entire function of exponential type $c$, 
and if $f$ and each of its derivatives have at 
least one zero in the unit circle, 
then $f$ is identically null.\\
In this paper we will use Gontcharoff estimate,
proved in \cite{Gon},
$$ 
\mid Q_n(x, x_0, x_1, \dots, x_{n-1})\mid \leq 
\frac{\left( \mid x - x_0\mid +\mid x_1 - x_0\mid +
\ldots + \mid x_{n-2} - x_{n-1}\mid\right)^n}{n!}.
$$
\subsection{Monotonicity property for some 
quasianalytic class}
\begin{proposition}
    Let $f \in C^\infty ([a, b])$ and a sequence
    $(x_n)_{n\in\N} $ such that 
    $ x_n\in [a,b], 
    \,\,\,\forall n\in\N$.
    Then, for each 
    $n \in \mathbb{N}$,
\begin{align*}
\begin{split}    
f(x) = f(x_0) + f'(x_1)Q_1(x, x_0) + \dots & + 
f^{(n)}(x_n)Q_n(x, x_0, \dots, x_{n-1})\\ &+ 
f^{(n+1)}(\xi)Q_{n+1}(x, x_0 , x_1, \dots, x_n),
\end{split}
\end{align*}
where $\xi$ is in the interval of ends 
$x_n$ and $x$.
\end{proposition}
\begin{proof}
    We put, for each $n \in \mathbb{N}$,
    \[
R_n(f)(x) = f(x) - f(x_0) - f'(x_1)Q_1(x, x_0) - 
\dots - f^{(n)}(x_n)Q_n(x, x_0, \dots, x_{n-1}). 
\]
By (6.5), for all $0 \leq k \leq n$, we have,
       \[
       (R_n(f))^{(k)}(x_k) = 
       f^{(k)}(x_k) - f^{(k)}(x_k)Q_k^{(k)}(x_k, x_0, x_1, \dots, x_{k-1}) = 0.
       \]
       and
       \[
       (R_n(f))^{(n+1)}(t) = f^{(n+1)}(t), \quad \forall t \in [a, b].
       \]
       We deduce then
       \[
       R_n(f)(x) = \int_{x_0}^x dt_1 \int_{t_1}^{x_1} dt_2 \dots \int_{t_{n-1}}^{x_{n-1}} dt_n \int_{t_n}^{x_n} (R_n(f))^{(n+1)}(t) dt.
       \]
       Hence there exists $\xi $  
       in the interval of ends $x_n$  and $x$, such that
       \[
       R_n(f)(x) = 
       f^{(n+1)}(\xi)Q_{n+1}(x, x_0 , x_1, \dots, x_n). 
       \]
    \end{proof}    
Let $n = (n_k)_k$ be a strictly increasing 
   sequence of natural numbers, and 
   $(x_n)_{n\in\N} $ a sequence such that 
    $ x_n\in [a,b], 
    \,\,\,\forall n\in\N$.

\begin{proposition}
Let $f\in C_{\overline{n}}([a, b])$ such that 
 $f^{(n)}(x_n) = 0$, for all $n \in \mathbb{N}$. 
If the series $\sum_{j=0}^\infty |x_j - x_{j+1}|$ 
converges, then the function 
$f$ is identically null.
\end{proposition}
\begin{proof}
    Since $f\in C_{\overline{n}}([a, b])$, 
    there exist two positive constants $A$, $B$, 
    such that
\[
M_{n_k}(f) \leq BA^{n_k} n_k!, \quad \forall k \in 
\mathbb{N}.
\]
Put $R_p = \sum\limits_{j=p}^\infty |x_j - x_{j+1}|$. 
Let $q_0 \in \mathbb{N}$, such that for all 
$q \geq q_0$, $R_q < \frac{1}{A}$. \\ We define 
the sequence $(m_s)_s$, by
\[
\forall s \in \mathbb{N}, \quad m_s + q + 1 \in 
\{n_k / k \in \mathbb{N}\}.
\]
It is clear that $m_s \rightarrow \infty$ 
when $s \rightarrow \infty$.
We apply Proposition 6.5 for $m_s + q$:
\[
f(x) = 
f^{(m_s+q+1)}(\xi)Q_{m_s + q+1}(x, x_0 , x_1, \dots, 
x_{m_s+q}).
\]
Taking $q$ times the derivative of each member, 
we obtain:
\begin{align*}
\begin{split}    
    f^{(q)}(x) &= f^{(m_s+q+1)}(\xi)
    Q^{(q)}_{m_s + q+1}(x, x_0, x_1, \dots, 
    x_{m_s +q})\\ 
    & = f^{(m_s+q+1)}(\xi)
    Q_{m_s +1}(x, x_q , x_{q+1}, 
    \dots, x_{m_s+q}).
\end{split}
\end{align*}
By using the Gontcharoff estimate for the 
polynomial $$  Q_{m_s +1}(x, x_q , x_{q+1}, 
\dots, x_{m_s+q}),$$ we have 
\begin{align*}
    \begin{split}    
|f^{(q)}(x)|\leq & B A^{m_s+q+1}(m_s+q+1)!\\&
\frac{\left(\mid x-x_q\mid + \mid x_{q+1}-x_q\mid 
+ \ldots + \mid x_{m_s + q- 1}-x_{m_s + q}\mid 
\right)^{( m_s + 1)}}
{( m_s + 1)!}\\ & = 
B A^{m_s+q+1}(m_s+q+1)!\frac{\left(\mid x-x_q\mid +
  R_q \right)^{ m_s + 1}}{( m_s + 1)!}\\ & =
  B A^{q}\frac{(m_s+q+1)!}{( m_s + 1)!}
  \left( A \mid x - x_q\mid + AR_q\right)^{m_s+1}.
\end{split}
\end{align*}
Let $\mu >0$, such that  
$\mu < \frac{1- AR_q}{A}$. We see that, 
for each $ x\in [x_q-\mu, x_q + \mu]$, we have 
$\left( A \mid x - x_q\mid + AR_q\right) < 1$. We 
deduce that the function 
$f^{(q)}$ is zero on the interval 
$[x_q - \mu, x_q + \mu]$, so 
 $f^{(q)}$ is zero 
on the interval $[a, b]$, by Proposition 6.1. 
If $q \geq  1$,  the function $f^{(q-1)}$ is constant 
on the interval $[a, b]$, hence zero, since 
$f^{(q-1)}(x_{q-1})=0$. By continuing with 
$f^{(q-2)},\ldots , f$,  we deduce that $f$ 
is   identically zero.
\end{proof}
\begin{remark}
We deduce that Corollary 4.4. is true for the 
quasianalytic class $C_{\overline{n}}([a, b])$,
and consequently Borel's Theorem 5.1.
on the non-surjectivity of the Borel 
mapping is also true.   
\end{remark}
Given the result of Theorem 5.2., it is natural 
to ask the following question:\\
{\bf{Problem}}\\
We know, by Theorem 5.2., that all quasianalytic 
classes do not satisfy Theorem 5.1 and consequently
corollary 4.4. We can then look for 
$C^\infty$ definable functions
in a polynomially bounded o-minimal structure.
More precisely, let   $\mathcal{R}$ be 
 a polynomially 
bounded o-minimal structure that extends 
the structure defined by the
restriction to the interval $[a, b]$ of analytic 
functions, see [7]. Is the monotonicity property 
verified by the $ C^\infty$ definable functions 
in this structure? What about Carleman's Theorem
about the non surjectivity of Borel mapping.
We have  an answer for the 
structure defined by the
multisommables functions, see [6, Corollary 8.6.].
\section{ appendix}
We say that a real function, $m$, of one real 
variable, is $C^\infty$ for $t\gg 0$, if there exists
$b >0$ such that the function $m$ is defined and 
$C^\infty$ in the interval $ [b,\infty[$. We write
$ m (t) >0$ for $t\gg 0$, if $ m (t) >0, \,\,\forall 
t\in [b,\infty[$. \\
In all the following, $m$ is a $C^\infty$ for 
$t\gg 0$. We suppose 
\begin{enumerate}
    \item $ m(t) > 0,\,\,\,m'(t) > 0,\,\,\, 
    m''(t) > 0,\,\,\,$ for $t\gg 0$,
    \item $\lim\limits_{t\to \infty}m'(t)= \infty$,
    \item there is $\delta >0$ such that 
    $ m''(t)\leq \delta $ for $t\gg 0$.
\end{enumerate}
We put $$ M(t) = e^{m(t)}\,\,\,\,\,\mbox{for}
\,\,\,\,\,t\gg 0.$$
\begin{definition}
A function $f$, $C^\infty$ in the interval $[a,b]$,
is said to be in the class $M$ in the interval 
$[a,b]$, if there exist $C >0$, $\rho >0$, such that
$$\mid f^{(n)}(t) \mid \leq C \rho^n M(n) \,\,\,\,
\mbox{for}\,\,\,\, n\gg 0, \,\,\,\,\forall t\in 
[a,b].$$
\end{definition}
We let $C_{M}([a,b])$  be the collection of all 
$C^\infty$ functions on $ [a,b]$ which are in 
the class $M$.
\begin{remark}
Let $M_1(t) = c r^t M(t)$, where $c>0,\,\, r >0$.
We easily see that a $C^\infty$ function, $f$, 
in the interval  $ [a,b]$ is in the class $M$, if 
and only if , $f$ is in the class $M_1$. Hence the 
class does not change when the function $ t\mapsto 
m(t)$ is changed by the function 
$ t\mapsto m(t) + \alpha t + \beta$, $ \alpha,\,
\beta\in \R$.
\end{remark}
We see then if $ m:[b,\infty[\to \R$, where 
$ b\geq 0$, we can suppose that $m(b)=0$, hence 
we extend $m$  to the interval $[0,\infty[$
by setting $m(t)=0$  if $t\leq  b$. We see that 
this extension is a convex function.
\begin{lemma}
 For all $j\in\N$, there exist $C_j >0,\,\, \rho_j >0$,
 such that 
 $$ M(p+j)\leq C_j\rho_j^p M(p),\,\,\,\,\,
 \forall p\in \N.$$  
\end{lemma}
\begin{proof}
 There exists $\theta\in ]p,p+j[$ such that   
 $m(p+j)- m(p) = jm'(\theta)$.
 Since $ m'' \leq \delta$, there exists $ C >0$ with
 $m'(t)\leq \delta t + C$. We have 
 $$ m'(\theta)\leq m'(p+j ) \leq \delta p + 
 ( C+\delta j).$$
 Put $\rho_j = e^{j\delta}$ and $ C_j = 
 e^{j(C +j\delta)}$, then we have 
 $M(p+j)\leq C_j\rho_j^p M(p)$.
\end{proof}
\begin{lemma}
    $C_{M}([a,b])$ is an algebra, closed under 
    differentiation.  
\end{lemma}
\begin{proof}
 Since the function $ m:[0,\infty [\to\R$ is convex 
 and  $m(0)=0$, we have, if $ 0\leq j\leq n$,
 $$ m(n-j)\leq \frac{n-j}{n}m(n),\,\,\,\,\mbox{and}
 \,\,\,\, m(j)\leq \frac{j}{n}m(n),$$
 hence $m(n-j) + m(j) \leq m(n)$, i.e.
 $M(n-j)M(j)\leq M(n)$. Using this inequality and 
 Leibnitz formula, we deduct the first statement 
 of the lemma.
 The second statement follows immediately 
 from lemma 7.3.
\end{proof}

 \begin{remark}
 If $ m(t) = t\log t$, i.e. $M(t) = t^t$, we have 
 the analytic class. If we want that the class 
 $C_M([a,b]) $ 
 cointains strictly the analytic class,   we choose
 the function $m$ of the form:
$$ m(t) = t\log t + t \mu(t),$$
where the function $\mu$ is increasing and 
$\lim\limits_{t\to\infty} \mu(t) = \infty$.\\
In order to have $m"(t) \leq \delta$, we must 
suppose that $\mu(t)\leq a t$, for $ t\gg 0$, 
where $a\in\R,\,\,\, a>0$.
We suppose also that $\mu$ is in a Hardy field.\\
Since the function $ t\mapsto t\mu(t)$ is convex,
the class $C_M([a,b]) $ is closed under 
composition, see \cite[2.5]{Child}.
\end{remark}
Let $t\mapsto M(t)= e^{m(t)}$ as above. For each 
$r >0$, we put 
$$ \Lambda(r) = \inf\limits_{ t\geq t_0}
\frac{M(t)}{r^t},$$
where $t_0$ is a fixed positive real.\\
The minimum of the function 
$t\mapsto \frac{M(t)}{r^t}$
is reached at a point $t$ verifying $ m'(t) =
\log r$ and this point is unique since the function
$t\mapsto m'(t)$ is increasing and 
$ \lim\limits_{t\to \infty}m'(t) = \infty$. We 
define a function $ r\mapsto \omega(r)$ by
$$\Lambda(r) = e^{-\omega(r)},
$$ we have 
$$\left\lbrace
\begin{array}{ll}
r= e^{m'(t)}\\
\omega(r)= tm'(t) + m(t),
\end{array}
\right.$$
or
$$
\left\lbrace
\begin{array}{ll}
r=et  e^{\mu(t)+t\mu'(t)}\\
\omega(r)= t + t^2\mu'(t).
\end{array}
\right.
$$
Since $\mu' >0$,  we see that $ \omega >0$,  and 
$\lim\limits_{t\to \infty}\omega(r) =\infty$.
We can easily inverse the last system; we have:
$$
\left\lbrace
\begin{array}{ll}
t =r \omega'(r)  \\
m(t)= r \omega'(r)\log r - \omega'(r).
\end{array}
\right.
$$
Since $m(t)= t\log t + t\mu(t)$, we have:
$$
\left\lbrace
\begin{array}{ll}
t =r \omega'(r)  \\
\mu(t)= -\log \omega'(r) -
\frac{\omega(r)}{ r \omega'(r)}.
\end{array}
\right.
$$
We see that the function $r\mapsto \omega(r)$
is increasing and when $ t\to \infty$, we have 
$$ r \omega'(r)\to\infty \,\,\,\,\,\,\mbox{and}
\,\,\,\,\,\,(\log \omega'(r) +
\frac{\omega(r)}{ r \omega'(r)})\to -\infty.
$$
We see then, when $r\to \infty$, we have 
$\omega'(r)\to 0$. For $ r >0$, we put
$$ \lambda(r)= 
\inf\limits_{n\in\N,\,\, n\geq  t_0}
\frac{M(n)}{r^n},$$
It is clear that $ \Lambda (r) \leq \lambda(r)$.
\begin{lemma}
    For $ r\gg 0$, we have 
$$e^{-\delta}\lambda(r)\leq \Lambda(r)\leq\lambda(r)
.$$
\end{lemma}
\begin{proof}
We put $$ \alpha(t) = m(t) - t\log r.$$
We see that 
$$ \Lambda(r) = e^{\alpha(t_0)},\,\,\,\,\mbox{with}
\,\,\,\,\alpha'(t_0) = 0.$$
Let  $n_0\in\N$ and $\mid n_0 - t_0\mid <1$, 
we have $\lambda(r) = e^{\alpha(n_0)}$. Note that
$$ \alpha(n_0) - \alpha(t_0) = 
\alpha'((1-\theta)n_0+\theta t_0) 
(n_0-t_0),\,\,\,\, 0< \theta < 1$$
and since $ m" \leq \delta$,
$$\mid \alpha'((1-\theta)n_0+\theta t_0) -
\alpha'(t_0)\mid \leq \delta,$$
we deduce then $e^{-\delta}\lambda(r)\leq 
\Lambda(r)$ and hence the proof of 
the lemma.
\end{proof}
\begin{proposition}
    The following three statements are equivalent:
\begin{enumerate}
\item $\sum\limits_{n=0}^\infty \frac{M(n)}{M(n+1)}
= \infty$
\item $\int_{r_0}^\infty \frac{\omega(r)}{r^2}dr= 
\infty$, for some $ r_0 >0$.
\item  $\int_{r_0}^\infty \frac{\log \lambda(r)}
{r^2}dr = -\infty$, for some $ r_0 >0$.
\end{enumerate}       
\end{proposition}
\begin{proof}
We have 
$$m'(n)\leq m(n+1) - m(n)\leq m'(n+1),$$
hence 
$$ \sum\limits_{n=0}^\infty \frac{M(n)}{M(n+1)} 
=\infty\Longleftrightarrow \int_{t_0}^\infty 
e^{-m'(t)}dt =\infty.
$$
Recall that by the above,
$$ \int_{t_0}^\infty e^{-m'(t)}dt = 
\int_{r_0}^\infty\frac{d(r\omega'(r))}{r}dr.$$
Since $\omega'(t)\to 0$ when $ r\to \infty$ and it 
is decreasing, we have 
$$\int_{r_0}^\infty\frac{d(r\omega'(r))}{r}dr 
= \infty \Longleftrightarrow
\int_{r_0}^\infty\frac{\omega(r)}{r^2}dr = 
\infty,$$ 
which proves $1. \Leftrightarrow  2.$\\
By Lemma 7.6., we have 
$$ -\frac{\log \lambda(r)}{r^2}\leq 
\frac{\omega(r)}{r^2}\leq \frac{\delta}{r^2} - 
\frac{\log \lambda(r)}{r^2},$$
hence $2. \Leftrightarrow  3.$
\end{proof}
By Theorem 2.4. and Theorem 2.8.,
if the class is quasianalytic, Proposition 7.7. 
tells us that the function $ r\mapsto \omega(r)$
tends to $\infty$ when $ r\to \infty$ as rapidly as
$r^q$, for all $ q<1$.{\it{ \underline {Probably the converse
of this statement is also true}}}.\\
In the case of the analytic class ($m(t)= t\log t$),
we have $ \omega(r)=  r \omega'(r)$, hence 
$\omega(r)= c r$, for some constant $c$. The 
converse is also true.
\begin{proposition}
If $ \omega (r)\backsimeq r$ when $r\to \infty$, then
any function $ f \in C_M([a,b])$ is analytic.   
\end{proposition}
\begin{proof}
By hypotheses, there exists $ c>0$ and $A >0$ such
that $$\omega(r) \geq cr,\,\,\, \forall r \geq A,$$
hence
$$ \forall p\in\N,\,\,\,  \forall r \geq A, \,\,\, 
e^{-\omega(r)}\leq \frac{p!}{(cr)^p}.$$
Since $m'(t)\to\infty$ when $ t\to \infty$, there
exists $N_0\in\N$ such that $ e^{m'(t)} \geq A$,
for all $ t\geq N_0$, ( we can suppose $N_0 >t_0$).
Let $ r > N_0$ and put $s= e^{m'(r)}$, we have then
$$ s\geq A\,\,\,\mbox{and}\,\,\,\frac{M(r)}{s^r}=
\Lambda(s)\leq \lambda(s)= \inf\limits_{n\in\N,\,\, n\geq  t_0}
\frac{M(n)}{r^n}.$$
By Lemma 7.6., we have $\lambda(s)\leq e^\delta 
\Lambda(s) = e^\delta e^{-\omega(s)}$, hence 
for all $r > N_0$
$$ \frac{M(r)}{s^r}\leq e^\delta e^{-\omega(s)}\leq
\frac{e^\delta}{ (cs)^p}p!,\,\,\,\forall p\in\N,$$
we obtain then 
$$ \forall p \in\N,\,\,\, p >N_0,\,\,\,
\frac{M(p)}{s^p}\leq \frac{e^\delta}{ (cs)^p}p!,$$
hence $$\forall p \in\N,\,\,\, p >N_0,\,\,\,
M(p)\leq \frac{e^\delta}{ c^p}p!.$$
 This proves the resul.
\end{proof}
\begin{proposition}
 Let $\mu(t)= \log\log t$, i.e. $m(t)= t\log t +
 t\log\log t$, then the class $C_M([a,b])$ is 
 quasianalytic, (recall $M(t) = e^{m(t)})$.
\end{proposition}
\begin{proof}
We will show that $\int_{r_0}^\infty 
\frac{\omega(r)}{r^2}dr= \infty$. \\
We have 
$$ s = e^{m'(t)} =e t\log t e^{\frac{1}{\log t}}
\backsimeq   e t\log t\,\,\,\mbox{when}\,\,\, t\to
\infty                      $$
and 
$$ \omega(s)= tm'(t) - m(t) = t + \frac{t}{\log t}
\backsimeq t \backsimeq \frac{s}{e\log s}.$$
Hence $$\frac{\omega(s)}{s^2} \backsimeq 
\frac{1}{es\log s},$$
which proves the proposition.
\end{proof}

\end{document}